\documentclass[leqno,12pt]{article}
\usepackage[margin=1in]{geometry}

\usepackage{dutchcal}

\usepackage{lmodern,microtype,mathtools,enumitem,needspace}

\usepackage{amsthm,amsfonts,amsbsy,amssymb,amsmath,amscd}

\usepackage[all]{xy}
\usepackage{xr-hyper}
\usepackage{tikz-cd}
\usepackage[colorlinks=true,linkcolor=blue,urlcolor=blue,citecolor=blue]{hyperref}

\usepackage[T1]{fontenc}

\usepackage{verbatim}
\usepackage{lscape}

\usepackage[textsize=footnotesize]{todonotes}

\newcommand{\BA}{{\mathbb {A}}}

\newcommand{\BC}{{\mathbb {C}}}

\newcommand{\BH}{{\mathbb {H}}}

\newcommand{\BN}{{\mathbb {N}}}

\newcommand{\BP}{{\mathbb {P}}}
\newcommand{\BQ}{{\mathbb {Q}}}
\newcommand{\BR}{{\mathbb {R}}}

\newcommand{\CH}{{\mathcal {H}}}

\newcommand{\CO}{{\mathcal {O}}}

\newcommand{\CZ}{{\mathcal {Z}}}

\newcommand{\an}{{\mathrm{an}}}

\newcommand{\Spec}{{\operatorname{Spec}}}

\newtheorem{thm}{Theorem}[subsection]
\newtheorem{maintheorem}{Theorem}

\newtheorem{cor}[thm]{Corollary}
\newtheorem{lem}[thm]{Lemma}
\newtheorem{prop}[thm]{Proposition}
\newtheorem {conj}[thm]{Conjecture}
\newtheorem{defn}[thm]{Definition}

\newtheorem{proposition}[thm]{Proposition}
\newtheorem{lemma}[thm]{Lemma}

\theoremstyle{definition}

\theoremstyle{remark}

\numberwithin{equation}{subsection}

\begin{document}
%------------------------------------------------------

\title{Heights and energies of polarized fibrations of varieties II\\[3pt]
\large Fibrations over a polarized base}
\author{Shou-Wu Zhang}
\date{\today}
\maketitle

\begin{abstract}
We extend the height and energy inequalities for fibrations over curves
to a fibration $f:X\to Y$ over a polarized variety $(Y,B)$.
For an integrable relatively nef adelic line bundle $\bar L$, we prove
a height inequality under the assumption that $L$ is nef on the inverse
image of a geometric generic complete-intersection curve in $Y$.
We deduce a lower bound for the local energy at almost every complex
or nonarchimedean realization of a fixed polarized model when $L$ is
$f$-ample. For a slope semistable bundle on a fixed polarized complex
variety, we show that the infimum of the energy on its projective bundle
equals the infimum of Donaldson's functional on Hermitian quotient
metrics, with the normalization given below.
\end{abstract}
\clearpage\begingroup\small\tableofcontents\endgroup\clearpage

\section*{Introduction}
\addcontentsline{toc}{section}{Introduction}

In Part~I we bounded heights and local energies for fibrations over a
curve. Let $K$ be a number field or the function field $k(T)$ of a smooth
geometrically integral projective curve over $k$. Let $f:X\to Y$ be a
flat projective morphism of smooth geometrically integral projective
$K$-varieties with geometrically integral generic fiber and
$\dim Y\ge2$. Fix nef adelic rational line bundles $\bar A_0$ on $X$
and $\bar B$ on $Y$ with ample underlying bundles. We put
\[
 e=\dim Y-1,\qquad d=\dim X-\dim Y+1,\qquad
 v_B(M)=M^d\cdot(f^*B)^e
\]
and take $\bar A=v_B(A_0)^{-1/d}\bar A_0$, so $A^d\cdot(f^*B)^e=1$.
We use relative nefness over $Y$ as in
\cite[Definition~1.3.3]{PartI}. Thus $\bar L$ is strongly relatively
nef if $\bar L+f^*\bar N$ is nef for some nef $\bar N$ on $Y$ with
$N$ ample. It is relatively nef if there is a nef adelic rational bundle
$\bar G$ on $X$ such that $\bar L+t\bar G$ is strongly relatively nef
for every rational $t>0$.
For very ample multiples $a_1B,\ldots,a_eB$, let $C$ be the
complete-intersection curve over the geometric generic point of
$\prod_i|a_iB|$, and write $X_C=X\times_Y C$.

\begin{maintheorem}\label{thm:main-height}
Let $\bar L$ be an integrable relatively nef adelic line bundle on $X$.
If $L|_{X_C}$ is nef for some choice of $a_1,\ldots,a_e$, then
\[
 \bar L^{d+1}\cdot(f^*\bar B)^e
 \le (d+1)v_B(L)(\bar L\cdot\bar A^d\cdot(f^*\bar B)^e).
\]
In particular, if $v_B(L)=0$, then $\bar L^{d+1}\cdot(f^*\bar B)^e\le0$.
\end{maintheorem}

For $f=\mathrm{id}_Y$ and $L\cdot B^e=0$, we recover the arithmetic
Hodge index inequality of Moriwaki and its adelic form due to Yuan
and Zhang \cite{Mo96,YZ17}. For $f:\BP_Y(E)\to Y$, where $E$ has
rank $r$, we apply Theorem~\ref{thm:main-height} to
$L=\CO(1)-r^{-1}f^*\det E$. Then $v_B(L)=0$. In characteristic
zero, if $E$ is geometrically slope semistable with respect to $B$,
its restriction to a general curve section of sufficiently high
multidegree is semistable \cite{MR82}, so $L|_{X_C}$ is nef.
In Section~\ref{ss:height-consequences} we obtain an adelic analogue
of Moriwaki's discriminant inequality \cite{Mo94,Mo95}.

We prove Theorem~\ref{thm:main-height} by induction on $\dim Y$. The
case of a curve is \cite[Theorem~A]{PartI}. For the induction step
we choose a divisor $H\in |mB|$ on which we can apply the theorem in
dimension $\dim Y-1$. Proposition~\ref{prop:mix-slicing} allows us to
recover intersection with $f^*\bar B$ by restricting to $f^{-1}(H)$
and dividing by $m$. Over a number field we use the theorem of Qu and
Yin \cite[Theorem~4.10]{QY24} to control the Green function of $H$.
Over $k(T)$ we choose $H$ on a projective model over $T$.

We denote the right side of Theorem~\ref{thm:main-height} minus the
left by $I_{A,B}(\bar L)$. To describe its variation, consider a tower
$X\xrightarrow{f}Y\to S$ of flat projective maps of pure relative
dimensions $d-1$ and $e+1$, and write $g:X\to S$ for the composite.
Put $v_{B,s}(L)=L_s^d\cdot(f_s^*B_s)^e$ and normalize $v_{B,s}(A)=1$.
We insert $e$ copies of $f^*B$ in each Deligne pairing of Part~I to
obtain the rational line on $S$
\[
 D_{A,B}(L)=
 -\langle L^{[d+1]},(f^*B)^{[e]}\rangle_g
 +(d+1)v_{B,s}(L)\langle L,A^{[d]},(f^*B)^{[e]}\rangle_g .
\]
Here $[j]$ denotes $j$ repeated entries. With the metrics on $A,B$
fixed, we write $D_{A,B}(h)$ for the metric induced by $h$ on $L$,
and $D_{A,B}(L,h)$ for the metrized line. When $S=\Spec K$ we have
$I_{A,B}(\bar L)=\widehat{\deg}\,D_{A,B}(\bar L)$. Over a valued
field, if $w_i$ is the weight of $D_{A,B}(h_i)$ in a common local
frame, we write $\Psi_{A,B}(h_0,h_1)=w_1-w_0$. Its value at $s$
is the local energy on $X_s$ given by \eqref{eq:H-energy}.

For the local bound let $F$ be a finitely generated extension of
$K_0=\BQ$ or $\mathbb F_p(t)$, and impose the same geometric
hypotheses on $f:X\to Y$ over $F$. In characteristic zero take
$E=\BC$ or a finite extension of $\BQ_p$, and in characteristic $p$
take $E$ to be a finite extension of $\mathbb F_p((t))$. At each
$K_0$-embedding $\sigma:F\hookrightarrow E$ we use the metrics on
$A,B$ from one polarized model (Definition~\ref{def:admissible-polarizations}).
We write $\Psi_\sigma=\Psi_{A_\sigma,B_\sigma}$ and denote by
$\BH^+(L_\sigma)$ the continuous relatively nef metrics.

\begin{maintheorem}\label{thm:main-energy}
Let $L$ be an $f$-ample rational line bundle such that $L|_{X_C}$ is
nef for a curve section as in Theorem~\ref{thm:main-height}.
For almost every $K_0$-embedding $\sigma:F\hookrightarrow E$ and every
integrable reference metric $h_{0,\sigma}$ on $L_\sigma$,
\[
 \inf_{h\in\BH^+(L_\sigma)}\Psi_\sigma(h_{0,\sigma},h)>-\infty .
\]
\end{maintheorem}

For a fixed polarized model,
one exceptional set works for all metrics in the theorem. The proof
uses the measures on the parameter space constructed in Part~I.
We show in Section~\ref{sec:families} how to apply Theorem~\ref{thm:main-height}
at algebraic parameters and control the extra factors of $f^*\bar B$
when we pass to their limiting measures.

For a fixed polarized complex variety $(Y,B)$ and a slope semistable
bundle $V$, Proposition~\ref{prop:semistable-projective-infimum}
gives a lower bound on $\BP_Y(V)$ itself. The infimum over all continuous
relatively nef metrics on the normalized tautological line equals the
infimum over Hermitian quotient metrics. On the latter the energy is
Donaldson's functional with the normalization of
\eqref{eq:normalized-bundle-donaldson}.

\paragraph*{Declaration on AI use.}
I used ChatGPT (OpenAI) to assist with language editing, exposition,
literature searches, and mathematical discussion. I also used it in
working out the cutoff and gluing arguments for energies in families
and the proof of Proposition~\ref{prop:semistable-projective-infimum},
where we use Donaldson's work to bound the energy on a semistable
projective bundle. I take full responsibility for the final text
and its mathematical content.

\phantomsection\label{par:notation-conventions}
\paragraph*{Conventions.}
We follow Part~I for local norms, Fubini--Study metrics, and adelic
degrees. In particular, $|z|_{\BC}=|z|^2$, and at a real place
$\Psi_{\BR}=\Psi_{\BC}/2$ for conjugation-invariant metrics.
We include the place weights in $\Psi_v$ and use the unnormalized
arithmetic degree over a number field. Over $k(T)$ we multiply
geometric intersection degrees by $c_k=\log\#k$ for finite $k$ and
by $c_k=1$ otherwise.

\section{The height inequality over a polarized base}\label{sec:heights}
We first prove Theorem~\ref{thm:main-height} by induction on $\dim Y$.
In Section~\ref{ss:height-consequences} we deduce the arithmetic
Hodge index inequality \cite{Mo96,YZ17} and an adelic discriminant
inequality. For projective bundles the latter gives the geometric
Bogomolov--Gieseker inequality.

\subsection{Induction on the dimension of the base}\label{ss:mix-polarized-base}

We argue by induction on $\dim Y$, allowing a geometrically integral
total space $X$ over a function field. This allows us to restrict $f$
to a smooth divisor of $Y$ even when its inverse image is singular.
Over a number field we retain the smoothness assumption on $X$.

For a divisor $H\in|mB|$ let $X_H=X\times_Y H$ and $f_H:X_H\to H$.
Write $B_H=B|_H$ and $U_H=U|_{X_H}$ for a line bundle $U$ on $X$.
Since $[X_H]=m f^*B$, the projection formula gives
\begin{equation}\label{eq:one-divisor-degree}
 v_{B_H}(U_H)=U_H^d\cdot(f_H^*B_H)^{e-1}=m\,v_B(U).
\end{equation}
In particular, the normalized polarization on $X_H$ is
$m^{-1/d}\bar A|_{X_H}$. Iterating the same formula for a curve $C$
cut out by divisors of degrees $k_1,\ldots,k_e$ in units of $B$ gives
\begin{equation}\label{eq:mix-degree}
 (U|_{X_C})^d=(k_1\cdots k_e)v_B(U).
\end{equation}
Relative nefness is preserved by base change and by integrable adelic
twists from the base \cite[discussion after Definition~1.3.3]{PartI}. To compare adelic
intersections on $X_H$ and $X$ over a number field, we first choose a
section defining $H$ with small Green function.

\subsubsection*{Choosing sections with small Green functions}

Let $\mathcal Z$ be a normal integral projective flat scheme over the
integers of a number field $F$, with smooth geometrically integral
generic fiber $Z$. Let $\overline{\mathcal N}$ be an arithmetically
ample smooth Hermitian line bundle on $\mathcal Z$, with generic fiber
$N$. For every sufficiently large integer $n$, fix a nonempty open
$U_n\subset |nN|$.

\begin{lem}\label{lem:mix-qu-yin}
For each $k\in\BN$, we can choose an integer $n_k\ge k$, a finite
extension $F_k/F$, and an integral section $s_k$ of $n_k\mathcal N$
over $\CO_{F_k}$ such that $[s_k]\in U_{n_k}(F_k)$ and $\operatorname{div}(s_k)$ has no
vertical components. For each embedding $\sigma:F\hookrightarrow\BC$,
\[
  \frac{1}{n_k}\log\|s_k^\tau\|\longrightarrow0
       \quad\text{in }L^1(Z_\sigma(\BC))\quad\text{as }k\to\infty,
\]
uniformly in the embeddings $\tau:F_k\hookrightarrow\BC$ extending
$\sigma$, where $\|\cdot\|$ denotes the pointwise Hermitian norm.
The same convergence holds after pullback by a fixed dominant morphism
$g:V\to Z_\sigma$ with $V$ smooth, projective, and integral.
\end{lem}

\begin{proof}
\smallskip\noindent\emph{Step 1. A family of sections meeting $U_{n_k}$.}
By \cite[Theorem~4.10]{QY24}, Qu and Yin give strictly increasing
integers $n_k$ and, for each $k$, an integral section $s$ of
$n_k\mathcal N$ over $\CO_F$. Their divisors are horizontal,
$\|s\|_{\rm sup}\to1$, and $n_k^{-1}\log\|s\|\to0$ in $L^1$.
We start with $n_1$ large enough that all the $U_{n_k}$ are defined. We perturb $s$ to place its divisor in
$U_{n_k}$ after a finite extension of $F$. By Zhang's bound for the last
successive minimum \cite[Theorem~3.10]{QY24},
we can complete $s$ to a basis $s,b_1,\ldots,b_h$ of $H^0(Z,n_kN)$
over $F$ with the $b_i$ integral and $\|b_i\|_{\rm sup}\leq\rho^{n_k}$
for some $0<\rho<1$ independent of $k$.

The basis $s,b_1,\ldots,b_h$ gives homogeneous coordinates
$[x_0:\cdots:x_h]$ on $|n_kN|$, where $h=h^0(Z,n_kN)-1$. On the
chart $x_0\ne0$, we divide by the coefficient $x_0$ of $s$ and
write the resulting section as $s+z_1b_1+\cdots+z_hb_h$, with
$z_i=x_i/x_0$. Since this chart is dense, it meets $U_{n_k}$ in a
nonempty open subset of $\mathbb A_F^h$. Every such open contains
the nonvanishing locus of some nonzero polynomial
$G\in F[z_1,\ldots,z_h]$. We may therefore choose $G$ so that
\[
 G(z_1,\ldots,z_h)\ne0
 \quad\Longrightarrow\quad [s+z_1b_1+\cdots+z_hb_h]\in U_{n_k} .
\]
We choose positive integers $p_i$ so that distinct monomials of
$G$ give distinct powers of $t$ after setting $z_i=t^{p_i}$. Then
$G(t^{p_1},\ldots,t^{p_h})$ is nonzero.

\smallskip\noindent\emph{Step 2. Specialization to a horizontal divisor.}
Consider $s(t)=s+\sum_i t^{p_i}b_i$. For fixed $k$, outside finitely
many primes of $F$ the special fiber of $\CZ$ is geometrically integral and
the basis remains independent after reduction and extension of the
residue field, by cohomology and base change. At the remaining
primes the constant term $s$ is nonzero on every geometric
component, since its divisor is horizontal. On each such component
only finitely many nonzero values of $t$ in the algebraic closure of
the residue field can make $s(t)$ identically zero. These values
have finite multiplicative order. Choose a primitive root $\zeta$ of
prime order greater than all those orders, prime to the finitely many
residue characteristics in question, and with
$G(\zeta^{p_1},\ldots,\zeta^{p_h})\ne0$. Reduction at these primes
preserves the order of $\zeta$. Taking $F_k=F(\zeta)$ and
$s_k=s(\zeta)$, we obtain a horizontal divisor with $[s_k]\in U_{n_k}(F_k)$.

\smallskip\noindent\emph{Step 3. $L^1$ convergence on $Z_\sigma$.}
Fix $\sigma:F\hookrightarrow\BC$ and write $c_{k,\sigma}=\|s^\sigma\|_{\rm sup}$.
The average of $n_k^{-1}\log\|s^\sigma\|$ against a smooth
probability measure tends to zero by $L^1$ convergence. It is at most
$n_k^{-1}\log c_{k,\sigma}$, which in turn is at most
$n_k^{-1}\log\|s\|_{\rm sup}\to0$. Thus
$n_k^{-1}\log c_{k,\sigma}\to0$. At every $\tau: F_k\hookrightarrow \BC$ above $\sigma$, the perturbation
has sup norm at most $h\rho^{n_k}$, where $h=O(n_k^{\dim Z})$. After division
by $c_{k,\sigma}$, this bound is still exponentially small. For any
sequence of embeddings $\tau_k:F_k\hookrightarrow\BC$ above $\sigma$,
we apply \cite[Lemma~3.9]{QY24} to $s^\sigma/c_{k,\sigma}$ and
$s_k^{\tau_k}/c_{k,\sigma}$. Adding back $n_k^{-1}\log c_{k,\sigma}$
gives $u_k:=n_k^{-1}\log\|s_k^{\tau_k}\|\to0$ in $L^1(Z_\sigma(\BC))$.
Since the sequence of embeddings was arbitrary, the convergence is uniform
over the embeddings above $\sigma$.

\smallskip\noindent\emph{Step 4. $L^1$ convergence after pullback to $V$.}
We now pull back $u_k$ by $g:V\to Z_\sigma$. Near any point of $V$ where
$dg$ is surjective, we choose holomorphic coordinates $(z,w)$ on $V$
and $z$ on $Z_\sigma$ in which $g(z,w)=z$. On a relatively compact
coordinate polydisc $U\times W$, Fubini's theorem gives
\[
 \int_{U\times W}|u_k\circ g|\,dV_z\,dV_w
 =\operatorname{vol}(W)\int_U |u_k(z)|\,dV_z \longrightarrow 0,
\]
where $dV_z$ and $dV_w$ denote Euclidean volume in these coordinates.
Thus $u_k\circ g\to0$ locally in $L^1$ on the submersion locus of $g$.

To pass to all of $V$, put $v_k=u_k\circ g$ and let $\omega$ be
the curvature form of the given Hermitian metric on $N_\sigma$.
With our normalization of $d^c$, the Poincar\'e--Lelong formula gives
\[
 2dd^c v_k
 =\frac{1}{n_k}\bigl[\operatorname{div}(g^*s_k^{\tau_k})\bigr]-g^*\omega
 \ \ge\ -g^*\omega ,
\]
since the divisor is effective. The lower bound is independent of $k$.
Fix one of the polydiscs $D=U\times W$ above, and normalize its
Euclidean volume measure $\mu_D$ to have mass one. Since
$v_k\le\sup_V v_k$ on $D$, we have
\[
 -\int_D |v_k|\,d\mu_D \le \int_D v_k\,d\mu_D
 \le \sup_V v_k \le \frac{1}{n_k}\log\bigl(\|s\|_{\rm sup}+h\rho^{n_k}\bigr).
\]
The left side tends to zero by $L^1$ convergence on $D$. The last
inequality follows from the triangle inequality for $s_k^{\tau_k}$,
and its right side tends to zero since $\|s\|_{\rm sup}\to1$ and
$h\rho^{n_k}\to0$. Hence $\sup_V v_k\to0$. By compactness for quasi-plurisubharmonic
functions \cite[Proposition~1.6]{GZ05}, every subsequence has a further
subsequence converging in $L^1(V)$. Its limit vanishes almost everywhere
on the submersion locus, whose complement has measure zero. Hence
$v_k\to0$ in $L^1(V)$.
\end{proof}

\subsubsection*{Restriction to a divisor}

Fix a rational $t>0$ and put $M=L+tA_0$. We choose a dense open
$Y^\circ\subset Y$ over which $f$ has geometrically integral fibers.
The integers $a_1,\ldots,a_e$ are those in Theorem~\ref{thm:main-height}.

\begin{lem}\label{lem:mix-replace-section}
For every positive multiple $m$ of $a_1$, there is a nonempty open
$U_m\subset |mB|$ of smooth geometrically integral divisors meeting
$Y^\circ$ with the following property. For $H\in U_m$, let $C_H$ be
the curve over the geometric generic point of
$\prod_{i=2}^e |a_iB_H|$, with $C_H=H$ if $e=1$. Then $M$ is ample
on $X_H\times_H C_H$. Over a number field we may also require $X_H$
to be smooth.
\end{lem}

\begin{proof}
On the inverse image of the geometric generic curve of multidegree
$(a_1,\ldots,a_e)$, the bundle $L$ is nef and $A_0$ is ample.
Thus $M$ is ample there. Ampleness on the fibers of the proper family
of inverse images of these curves is open \cite[Tag 0D2S]{Stacks}.
By Bertini and the infinitude of $K$, we may therefore choose divisors
$D_i\in |a_iB|$ over $K$ such that $W=D_2\cap\cdots\cap D_e$ is a
smooth geometrically integral surface, $C_0=D_1\cap W$ is a smooth
geometrically integral curve meeting $Y^\circ$, and $M|_{X_{C_0}}$
is ample. For $e=1$ we take $W=Y$.

Write $m=qa_1$. On $W$, the divisor $qD_1$ cuts out a thickening of
$C_0$, and its inverse image in $X$ is a nilpotent thickening of
$X_{C_0}$. Ampleness is unchanged by such a thickening. Applying
openness of ampleness to the family $X_{W\cap H}$ for $H\in |mB|$,
we obtain an open neighborhood of $qD_1$ on which $M|_{X_{W\cap H}}$
is ample. We intersect this neighborhood with the Bertini open where
$H$ and $W\cap H$ are smooth and geometrically integral and $W\cap H$ meets $Y^\circ$.
In characteristic zero, Bertini applied to the base point free system
$f^*|mB|$ on a smooth $X$ also gives smoothness of $X_H$.

For any $H$ in the resulting open denoted by $U_m$, the divisors $D_2|_H,\ldots,D_e|_H$
cut out $W\cap H$. If $e>1$, openness of ampleness in the product
of the full systems $|a_iB_H|$ therefore gives ampleness on the
inverse image of $C_H$. For $e=1$ we already have ampleness on $X_H$.
\end{proof}

For $H\in U_m$, the morphism $f_H$ is flat by base change and has
geometrically integral generic fiber since $H$ meets $Y^\circ$.
It follows by flatness that $X_H$ is geometrically integral.

We compare intersections over $K$ even when $H$ is defined over a
finite extension $F$. For a number field $K$, write $I_{X_H/K}$ for
the adelic intersection on $X_H$ over $F$ divided by $[F:K]$. For
$K=k(T)$ we use only constant extensions $F=k'(T)$ and multiply the
adelic intersection over $F$ by $c_k/c_{k'}$. On a model this is the
geometric intersection degree over $k'$ multiplied by $c_k$.

\begin{prop}\label{prop:mix-slicing}
Put $D=\dim X$. Let $Q=\prod_{i=1}^D\bar Q_i$ run over a finite
set of products of integrable adelic line bundles on $X$. For every $\delta>0$ there are a multiple
$m$ of $a_1$, a finite extension $F/K$ as above, and a divisor
$H\in U_m(F)$ such that
\[
 \left|\frac{1}{m}I_{X_H/K}(Q|_{X_H})-Q\cdot f^*\bar B\right|<\delta
 \qquad\text{for every }Q .
\]
\end{prop}

\begin{proof}
We first work on a flat projective model $\mathcal X\to\mathcal Y$
of $f$ with $\mathcal Y$ normal. Assume that the factors of $Q$ have
model metrics and that $\bar B$ comes from an ample rational bundle
$\mathcal B$ on $\mathcal Y$. Over a number field we assume that
$\overline{\mathcal B}$ is arithmetically ample with smooth positive
curvature and that the metrics in $Q$ are smooth at infinity.

\smallskip\noindent\emph{Number fields.}
We choose a multiple $r$ of $a_1$ such that $r\mathcal B$ is a line
bundle and apply Lemma~\ref{lem:mix-qu-yin} on $\mathcal Y$ to
$r\overline{\mathcal B}$ and the open subsets $U_{nr}\subset |nrB|$.
Let $s$ be a section supplied by the lemma over $F$, put $m=nr$,
and write $H$ for its divisor on $Y_F$. Its model divisor has no
vertical components. The same holds after flat pullback to the model
of $X_F$.

For $Q=\bar Q_1\cdots\bar Q_{\dim X}$ and
$\sigma:K\hookrightarrow\BC$, put
$\Omega_{Q,\sigma}=\prod_i c_1(\bar Q_{i,\sigma})$ on $X_\sigma(\BC)$.
The arithmetic intersection formula gives
\[
 Q\cdot f^*\bar B-\frac{1}{m}I_{X_H/K}(Q|_{X_H})
 =-\frac{1}{[F:K]}\sum_{\tau:F\hookrightarrow\BC}
   \int_{X_\sigma(\BC)}
   \frac{\log\|f_\sigma^*s^\tau\|}{m}\,\Omega_{Q,\sigma},
 \qquad \sigma=\tau|_K .
\]
Here $\|\cdot\|$ is the ordinary Hermitian norm.
Lemma~\ref{lem:mix-qu-yin} applied to $f_\sigma$ gives $L^1$
convergence to zero of the normalized logarithms, uniformly over the
$[F:K]$ embeddings above each $\sigma$. Since the forms
$\Omega_{Q,\sigma}$ are fixed, we may choose one section giving error
less than $\delta$ for every $Q$.

\smallskip\noindent\emph{Function fields.}
Let $p:\mathcal Y\to T$ be the structure morphism. We choose a
multiple $r$ of $a_1$ such that $r\mathcal B$ gives a closed immersion
$i:\mathcal Y\hookrightarrow\BP^N_k$. Since $T$ is separated over $k$,
the graph of $p$ is closed. Thus $(i,p)$ embeds $\mathcal Y$ as a
closed subscheme of $\BP^N_k\times_k T$. Restricting to the generic
point of $T$ gives a closed immersion $j:Y\hookrightarrow\BP^N_K$
with $j^*\mathcal O(1)=rB$.

Let $\mathcal I_Y$ be the ideal sheaf of $j(Y)$ in $\BP^N_K$. For
$n$ sufficiently large, Serre vanishing gives
$H^1(\BP^N_K,\mathcal I_Y(n))=0$. The ideal sequence twisted by
$\mathcal O(n)$ then gives a surjection
\[
 K[x_0,\ldots,x_N]_n=H^0(\BP^N_K,\mathcal O(n))
 \xrightarrow{\,j^*\,} H^0(Y,nrB).
\]
Write $V_n=H^0(\mathcal Y,nr\mathcal B)$. Each monomial
$x_0^{u_0}\cdots x_N^{u_N}$ with $\sum u_i=n$ is defined over $k$,
so its restriction by $i$ belongs to $V_n$. The images of these
sections on $Y$ span $H^0(Y,nrB)$ over $K$ by the displayed
surjection. Hence the restriction map
$\rho_n:V_n\otimes_k K\to H^0(Y,nrB)$ is surjective.

The model $\mathcal Y$ is geometrically integral over $k$ by flatness
over $T$ and geometric integrality of $T$ and $Y$. Choose a $k$-basis
$b_1,\ldots,b_h$ of $V_n$ and write $s=\sum_{i=1}^h z_i b_i$.
By Bertini the vectors $z$ for which $\operatorname{div}(s)$ is
geometrically integral form a nonempty open in $\BA^h_k$.
Over $K$ the conditions $\rho_n(s)\ne0$ and $[\rho_n(s)]\in U_{nr}$
also define a nonempty open in $\BA^h_K$, since $\rho_n$ is surjective.
The two opens meet after base change to $K$, since $\BA^h_K$ is
irreducible. Their intersection contains $\{G\ne0\}$ for some
nonzero $G\in K[z_1,\ldots,z_h]$.

We take $k'=k$ if $k$ is infinite. If $k$ is finite, we choose a
finite extension $k'/k$ with $\#k'>\deg G$. In either case there is
a tuple $(c_1,\ldots,c_h)\in(k')^h$ with $G(c_1,\ldots,c_h)\ne0$
in $k'(T)$. Indeed, $G$ remains nonzero over $k'(T)$, and a polynomial
of degree less than $\#k'$ in each variable that vanishes on $(k')^h$
is zero, by the one-variable root bound applied successively.
Thus $s=\sum_{i=1}^h c_i b_i$ has geometrically integral model divisor
and its generic fiber $H$ belongs to $U_{nr}$. In particular, the
model divisor is horizontal. Its class is $m\mathcal B$ for $m=nr$, so flat pullback
and the projection formula give the equality
\[
 \frac{1}{m}I_{X_H/K}(Q|_{X_H})=Q\cdot f^*\bar B .
\]

\smallskip\noindent\emph{Passage to adelic metrics.}
By multilinearity we may work with strongly nef factors of $Q$ with
ample underlying bundles. Indeed, each integrable factor is a
difference of strongly nef bundles, and we can add the same ample
model bundle with nef metric to both terms. We approximate the factors
by nef models on the same underlying bundles. At infinity we use
smooth semipositive metrics as in the proof of
\cite[Theorem~1.3.4]{PartI}.

We also approximate $\bar B$ on its underlying bundle $B$.
By nefness there is a strongly nef $\bar G$ on $Y$ such that
$j\bar B+\bar G$ is strongly nef for every $j\ge1$
\cite[\S2.5]{YZ26}. For some integer $b$ the bundle $bB-G$ is ample.
We give $bB-G$ a nef model metric and denote the resulting bundle by $\bar E$.
Then $(j\bar B+\bar G+\bar E)/(j+b)$ is strongly nef on $B$ and
converges to $\bar B$. Thus we can first approximate $\bar B$ by nef
models on $B$.

To obtain the ample models used above, we start with a projective
model of $Y$ carrying an ample extension of a multiple of $B$. On a
common model with the nef approximation to $\bar B$, we take the
closure of the graph of $f$ using a model carrying the factors of $Q$.
After a blowup of the base away from the generic fiber the strict
transform of the graph is flat. A blowup of the initial polarized
model dominates this base modification \cite[Tags 081R and 081T]{Stacks}.
Let $\pi$ denote this blowup, $E$ its exceptional divisor, and
$\mathcal B_0$ the original ample extension of $b_0B$. For $l$ large,
$\pi^*(l\mathcal B_0)(-E)$ is ample. Dividing by $lb_0$ gives an
ample rational model of $B$. We may normalize the base and
pull back the flat morphism. Over a number field we choose a smooth
positive metric on the ample model and scale the norms to make it
arithmetically ample. We thus obtain the required ample approximation to $\bar B$
by taking a sufficiently small positive proportion
of this model in a convex combination with the nef approximation.

We now bound the error uniformly in $H$. Suppose we replace the
metric on $Q_1$ by a second metric on the same bundle, keeping the
nef factors $\bar Q_2,\ldots,\bar Q_D$ fixed. At each place $v$ of $K$
let $\phi_v=\log(|\cdot|_{\rm old,v}/|\cdot|_{\rm new,v})$ on
$X_v^{\rm an}$, and put $\eta=\sum_v\|\phi_v\|_{\sup}$. Here the
absolute values already include the place weights.

Let $H\in |mB|$ be defined over $F$. For each place $w$ of $F$
above $v$, we pull $\phi_v$ back to $X_{H,w}^{\rm an}$ and use the
mixed measure
\[
 \mu_{H,w}=c_1(\bar Q_{2,H})_w\cdots c_1(\bar Q_{D,H})_w .
\]
The measure is positive since the metrics are semipositive. Its mass
is the geometric intersection number on $X_H$
\cite[Appendix~1]{YZ17}. As $[X_H]=m f^*B$ on $X_F$, we have
\[
 \int_{X_{H,w}^{\rm an}}\mu_{H,w}
  =(Q_{2,H}\cdots Q_{D,H})=m(Q_2\cdots Q_D\cdot f^*B).
\]
With the normalization of $I_{X_H/K}$ above, the change-of-metric
formula gives the contribution at $v$ as
\[
 \sum_{w\mid v}\frac{[F_w:K_v]}{[F:K]}
       \int_{X_{H,w}^{\rm an}}\phi_v\,d\mu_{H,w}.
\]
The coefficients sum to one. Each integral is bounded in absolute
value by $\|\phi_v\|_{\sup}$ times the mass of $\mu_{H,w}$.
Dividing by $m$ and summing over $v$, we see that the change in
$m^{-1}I_{X_H/K}(Q|_{X_H})$ is at most
\[
 \sum_v\|\phi_v\|_{\sup}\,(Q_2\cdots Q_D\cdot f^*B)
   =\eta\,(Q_2\cdots Q_D\cdot f^*B)
\]
in absolute value, independently of $H$, $m$, and $F$.
On a projective variety
the weighted sup errors tend to zero under model approximation
\cite[\S2.5 and Appendix~A.5]{YZ26}. Changing one metric at a time
therefore gives error less than $\delta/3$ in the normalized
intersection on $X_H$, uniformly in $H$. By continuity we may also
make the error in $Q\cdot f^*\bar B$ less than $\delta/3$ for every $Q$. For these models we choose $H$
with error less than $\delta/3$, which proves the proposition.
\end{proof}

\begin{proof}[Proof of Theorem~\ref{thm:main-height}]
When $\dim Y=1$, the hypothesis says that $L$ is nef on $X$, and the
inequality is \cite[Theorem~A]{PartI}. Suppose $\dim Y\ge2$ and
assume the theorem in dimension $\dim Y-1$. Fix a rational $t>0$
and put $\bar M=\bar L+t\bar A_0$. To pass from a divisor $H$ to
$Y$, we apply Proposition~\ref{prop:mix-slicing} to the two products
\[
 Q_1=\bar M^{d+1}(f^*\bar B)^{e-1},\qquad
 Q_2=\bar M\bar A^d(f^*\bar B)^{e-1}.
\]
Here $\bar A^d=\bar A_0^d/v_B(A_0)$. For each $\delta>0$ the
proposition gives $H\in U_m(F)$ approximating both intersections.
By Lemma~\ref{lem:mix-replace-section}, $M$ is ample on the inverse
image of the geometric generic curve on $H$. The restriction
$\bar M_H$ is relatively nef over $H$. We can therefore apply the
induction hypothesis to $f_H:X_H\to H$ with base polarization
$\bar B_H$ and normalized polarization $m^{-1/d}\bar A|_{X_H}$.

By \eqref{eq:one-divisor-degree}, $v_{B_H}(M_H)=m v_B(M)$. The
factor $m$ cancels the factor $m^{-1}$ from the $d$ copies of the
normalized polarization on $X_H$. Dividing the inequality over $H$
by $m$ and expressing both sides over $K$, we obtain
\[
 \frac{1}{m}I_{X_H/K}(Q_1|_{X_H})
 \le (d+1)v_B(M)\,\frac{1}{m}I_{X_H/K}(Q_2|_{X_H}).
\]
Letting $\delta\to0$ gives
$Q_1\cdot f^*\bar B\le(d+1)v_B(M)(Q_2\cdot f^*\bar B)$.
This is the desired inequality for $\bar M$. Both sides are
polynomial in $t$, so letting $t\downarrow0$ proves the theorem.
\end{proof}

Iterating Proposition~\ref{prop:mix-slicing} also gives approximation
on a curve. Indeed, after each restriction $M$ is ample on the inverse
image of the geometric generic curve of the new base, which is the
ampleness condition used in Lemma~\ref{lem:mix-replace-section}.
For a finite set of products $R$ of $d+1$ integrable factors, we
apply the proposition at the $j$th step to the restrictions of
$R(f^*\bar B)^{e-j}$, with error less than $\delta/e$ after
normalization over $K$. If the successive divisors have degrees $m_1,\ldots,m_e$ in units of $B$ and
their final intersection is $C$, the $j$th error enters with weight
$(m_1\cdots m_{j-1})^{-1}\le1$. Hence
\[
 \left|\frac{I_{X_C/K}(R|_{X_C})}{m_1\cdots m_e}
                 -R\cdot(f^*\bar B)^e\right|<\delta .
\]

\subsection{The Hodge index and discriminant inequalities}\label{ss:height-consequences}

Let $X$ be a smooth projective geometrically integral variety over a
number field, and take $f=\mathrm{id}_X$. Then $d=1$ and $e=\dim X-1$.
Every integrable adelic line bundle $\bar L$ is strongly
relatively nef for the identity. Indeed, write $\bar L=\bar M-\bar N$
with $\bar M,\bar N$ nef. Adding the same nef adelic polarization to
$\bar M$ and $\bar N$, we may assume that $N$ is ample. If $L\cdot B^e=0$, then
$L|_C$ has degree zero by \eqref{eq:mix-degree}, so it is nef.
Theorem~\ref{thm:main-height} gives
\[
 \bar L^2\cdot\bar B^e\le 0.
\]
For Hermitian line bundles on arithmetic models, this is the inequality
in Moriwaki's arithmetic Hodge index theorem \cite[Theorem~B]{Mo96}.
For adelic line bundles it follows from the Hodge index theorem of
Yuan and Zhang \cite[Theorem~3.2]{YZ17} by taking all the polarizations
equal to $\bar B$.
The curve case goes back to Faltings and Hriljac \cite{Fa84,Hr85}.

Let $f:X\to Y$ satisfy the geometric hypotheses of
Theorem~\ref{thm:main-height}, and let $\eta$ be the generic point
of $Y$. For an integrable adelic bundle $\bar L$ on $X$ with
$L_\eta^{d-1}>0$, the Deligne pairing
$D_f(\bar L)=\langle\bar L^{[d]}\rangle_{X/Y}$ is an integrable
adelic bundle on $Y$ \cite[Theorem~4.1.3]{YZ26}. As in Part~I,
we define the normalized bundle and the relative discriminant by
\begin{align*}
 \bar L^\circ&=\bar L-\frac{f^*D_f(\bar L)}{dL_\eta^{d-1}},\\
 \Delta_{f,\bar B}(\bar L)
  &=\frac{d+1}{2dL_\eta^{d-1}}\bigl(D_f(\bar L)^2\cdot\bar B^e\bigr)
       -\bar L^{d+1}\cdot(f^*\bar B)^e .
\end{align*}

\begin{thm}\label{thm:bgi}
Let $\bar L$ be relatively nef over $Y$. If $L^\circ$ is nef on
$X_C$ for a curve section as in Theorem~\ref{thm:main-height}, then
$\Delta_{f,\bar B}(\bar L)\ge0$.
\end{thm}
\begin{proof}
As in the proof of \cite[Theorem~1.4.1]{PartI}, $\bar L^\circ$ is invariant under
base twists and is relatively nef. In higher dimension the numerical
identity we need on $Y$ is $f_*((L^\circ)^d)=0$. Indeed,
$f_*(L^d)=c_1(D_f(L))$ and $f_*(L^{d-1})=L_\eta^{d-1}[Y]$.
Intersecting with
$B^e$ gives $v_B(L^\circ)=0$. With $e$ copies of $f^*\bar B$,
the same projection formula gives
\[
 (\bar L^\circ)^{d+1}\cdot(f^*\bar B)^e
 =\bar L^{d+1}\cdot(f^*\bar B)^e
  -\frac{d+1}{2dL_\eta^{d-1}}\bigl(D_f(\bar L)^2\cdot\bar B^e\bigr)
 =-\Delta_{f,\bar B}(\bar L).
\]
Theorem~\ref{thm:main-height} applied to $\bar L^\circ$ proves the
inequality and the identity also gives invariance of the
discriminant under base twists.

\end{proof}

Moriwaki proved the arithmetic Bogomolov--Gieseker inequality for
Hermitian vector bundles on arithmetic surfaces \cite{Mo94} and on
higher-dimensional polarized arithmetic varieties \cite{Mo95}. In the
latter case he assumes slope semistability on each connected component
of the complex fiber. We now apply Theorem~\ref{thm:main-height} to a
projective bundle over $(Y,B)$.

For a vector bundle $E$ of rank $r\ge2$ on a polarized base $Y$,
let $f:\BP_Y(E)\to Y$ parametrize its quotient lines, and put
$\xi=\CO(1)$ and $L=\xi-r^{-1}f^*\det E$. Here $d=r$.
The projective bundle formula gives $f_*\xi^{r-1}=1$ and
$f_*\xi^r=c_1(E)$, so $f_*L^r=0$.

\begin{cor}\label{cor:mix-projective-bundle}
Suppose that on a geometric generic curve section $C$ the restriction
$E|_C$ remains semistable after every finite cover of smooth projective
curves. Then every integrable relatively nef adelic metric on $L$
satisfies
\[
 \bar L^{r+1}\cdot(f^*\bar B)^e\le0 .
\]
In characteristic zero it suffices that $E$ be geometrically slope
semistable with respect to $B$.
\end{cor}
\begin{proof}
A curve in $\BP_C(E|_C)$ dominating $C$ gives on its normalization
a finite map $g:C'\to C$ and a quotient line $Q$ of $g^*(E|_C)$.
Semistability gives $\deg Q\ge(\deg g)\deg(E|_C)/r$. This is precisely
the nef inequality for $L$ on that curve. On a curve in a fiber,
$L$ has positive degree. Since $v_B(L)=0$, the assertion follows
from Theorem~\ref{thm:main-height}. In characteristic zero we use the
restriction theorem of Mehta and Ramanathan \cite{MR82} and preservation
of semistability under finite pullback on a smooth curve. In positive
characteristic the hypothesis includes pullback by Frobenius.
\end{proof}

For $e\ge1$, we can read a geometric consequence from this inequality.
Choose a strongly relatively nef metric on $L$ by writing it as
$(L+f^*N)-f^*N$ with $L+f^*N$ and $N$ ample and equipping both
ample bundles with nef metrics. Let $\bar U$ be a line from $\Spec K$
of positive adelic degree. Replacing $\bar B$ by $\bar B+s\bar U$
for rational $s\ge0$, we obtain an affine function of $s$ which
is nonpositive. Its linear coefficient is
$e\,\widehat{\deg}(\bar U)\,(L^{r+1}\cdot (f^*B)^{e-1})$.
Since $f_*\xi^{r+1}=c_1(E)^2-c_2(E)$, we obtain
\[
 f_*L^{r+1}=\frac{r-1}{2r}c_1(E)^2-c_2(E),\qquad
 \big(2r c_2(E)-(r-1)c_1(E)^2\big)\cdot B^{e-1}\ge0.
\]

\section{Energies of fibrations}\label{sec:energy}
We define the local energy over a valuation field $F$. When the
fibration is defined over a global field, we use Theorem~\ref{thm:main-height}
to bound the energy below at each place.

\subsection{Local metrics and the mixed energy}\label{ss:local-metrics}

Let $f:X\to Y$ be a flat projective fibration of smooth projective
varieties over a valuation field $F$. We take analytifications over its
completion and equip $A,B$ with semipositive metrics such that $v_B(A)=1$. Put $\omega_A=c_1(\bar A)$,
$\beta=c_1(f^*\bar B)$, and $\alpha_i=c_1(L,h_i)$ for continuous
integrable metrics $h_0,h_1$ on $L$. The change-of-metric formula for
$D_{A,B}(L)$ gives
\begin{equation}\label{eq:H-energy}
 \Psi_{A,B}(h_0,h_1)=
 \int_{X^{\rm an}}\log\frac{|\cdot|_{h_1}}{|\cdot|_{h_0}}
 \left(\sum_{i=0}^d\alpha_0^i\alpha_1^{d-i}
       -(d+1)v_B(L)\omega_A^d\right)\beta^e .
\end{equation}
The measure has mass zero, so rescaling either metric by a constant
preserves the energy. As a difference of weights, $\Psi_{A,B}$ satisfies
the cocycle identity. For $v_B(L)=0$ it is independent of $\bar A$.

We use the metric spaces of \cite[\S2.1]{PartI}. A metric is strongly
relatively nef if there is an ample bundle $N$ on $Y$ with a
continuous semipositive metric $n$ such that
\begin{equation}\label{eq:admissible-base-twist}
 (L,h)+f^*(N,n)\ \text{is semipositive on }X.
\end{equation}
It is relatively nef if it is continuous and integrable and adding
$\epsilon\bar G$ makes it strongly relatively nef for every rational
$\epsilon>0$, for one fixed semipositive rational bundle $\bar G$ on $X$.
We write $\BH^+(L)$ for these metrics. For $f$-ample $L$, $\CH^+(L)$
denotes the smooth metrics strictly positive on $\ker df$ over $\BC$,
and the rational model metrics from relatively ample bundles over a
nonarchimedean field. An ample base twist shows that
$\CH^+(L)\subset\BH^+(L)$, also at singular fibers, as in
\cite[\S2.1]{PartI}. The same construction gives a strongly relatively
nef metric on an $f$-semiample bundle by writing a positive multiple
as the difference of a globally generated bundle and an ample base twist.

\begin{lem}\label{lem:relative-nef-approximation}\label{lem:relative-nef-regularization}
Suppose that $L$ is $f$-ample and $Y$ is smooth. Over $\BC$ suppose
also that $X$ is smooth. Over any complete nontrivially valued
nonarchimedean field we allow singular $X$. Each $h\in\BH^+(L)$ is a uniform limit of
metrics $h_j\in\CH^+(L)$ on $L$ such that, for every integrable $h_0$,
\[
 \Psi_{A,B}(h_0,h_j)\longrightarrow\Psi_{A,B}(h_0,h).
\]
For discrete valuations the approximating models can be chosen projective
and flat over projective models of $Y$.
\end{lem}
\begin{proof}
The proofs of \cite[Lemmas~2.1.1 and~2.1.2]{PartI} use the base
only to twist by an ample bundle and to form the graph of $f$ on
models. They therefore give the same uniform approximants over $Y$.
For a strongly relatively nef metric the approximants are differences
of semipositive metrics on two fixed bundles. In their mixed
intersections we keep the $e$ copies of $f^*\bar B$ fixed, so
continuity under uniform convergence gives the asserted energy limit
\cite[\S4.2.2]{YZ26}, \cite[Theorem~8.16]{BE21}.
For a relatively nef metric the proof of \cite[Lemma~2.1.2]{PartI}
first interpolates on $L$ between two fixed integrable metrics.
The mixed energy along this segment is polynomial by multilinearity.
The same diagonal choice then gives convergence for one reference
metric, hence for every reference metric by the cocycle identity.
\end{proof}

In particular the infima on $\CH^+(L)$ and $\BH^+(L)$ agree.
For smooth $f$ over $\BC$, interpolation with a fixed strictly
vertically positive metric also gives the same infimum on the smooth
fiberwise-semipositive metrics. We formulate the mixed version of
\cite[Conjecture~2.1.3]{PartI} for $F=\BC$ or a complete discretely
valued field. We impose the geometric hypotheses of
Theorem~\ref{thm:main-height} over $F$, allowing also a curve base
with $C=Y$ as in Part~I.

\begin{conj}\label{conj:H-energy-boundedness}
Let $L$ be an $f$-ample rational line bundle whose restriction to
$X_C$ is nef for some geometric generic curve section $C$.
Then $\inf_{h\in\CH^+(L)}\Psi_{A,B}(h_0,h)>-\infty$ for every
integrable reference metric $h_0$.
\end{conj}

For a smooth path of metrics over $\BC$ put
$\psi_t=\log(|\cdot|_{h_t}/|\cdot|_{h_0})$ and $\alpha_t=c_1(L,h_t)$.
In the variation formulas of \cite[\S2.1 and (2.4.1)]{PartI}, the
additional closed factor $\beta^e$ is fixed under integration by parts.
Writing $s=v_B(L)$, we obtain
\begin{equation}\label{eq:H-energy-variations}
 \begin{aligned}
 \frac{d}{dt}\Psi_{A,B}(h_0,h_t)
    &=(d+1)\int_X\dot\psi_t(\alpha_t^d-s\omega_A^d)\beta^e,\\
 \frac{d^2}{dt^2}\Psi_{A,B}(h_0,h_t)
    &=-d(d+1)\int_X\psi\,dd^c\psi\,\alpha_t^{d-1}\beta^e
       \qquad(\psi_t=t\psi).
 \end{aligned}
\end{equation}
Along an integrable affine variation $\psi_t=t\psi$ the same formulas
follow by differentiating the mixed intersection polynomial, also for
nonarchimedean model metrics. Strong relative nefness along the
segment follows by taking a common semipositive base twist at the
endpoints. Testing the first variation against smooth functions or
model functions, respectively, gives the critical equation
\begin{equation}\label{eq:mixed-critical-equation}
 \alpha^d\beta^e=v_B(L)\omega_A^d\beta^e .
\end{equation}

Suppose that $L$ is $f$-ample, $v_B(L)=0$, and $L^\circ|_{X_C}$
is nef for a geometric generic curve section $C$ as in
Theorem~\ref{thm:main-height}. Every globally semipositive metric $h$
on $L$ satisfies \eqref{eq:mixed-critical-equation}, since the positive
measure $c_1(L,h)^d\beta^e$ has total mass $v_B(L)=0$.
For a fixed integrable reference metric $h_0$, which solutions
$h_*\in\BH^+(L)$ of \eqref{eq:mixed-critical-equation} minimize
$\Psi_{A,B}(h_0,\cdot)$ on $\BH^+(L)$? For a curve base, we formulated
the corresponding conjecture in \cite[Conjecture~2.4.3]{PartI}.

\begin{proposition}\label{prop:zero-volume-energy}
Let $\eta$ be the generic point of $Y$ and suppose that
$v_B(L)=L_\eta^{d-1}=0$. Over $\BC$ or a complete nontrivially
valued nonarchimedean field, two continuous strongly relatively nef
metrics on $L$ satisfy
\[
 \alpha_0^i\alpha_1^{d-i}\beta^e=0\quad(0\le i\le d),
 \qquad\Psi_{A,B}(h_0,h_1)=0 .
\]
For a smooth complex fibration the same holds for two smooth
fiberwise-semipositive metrics.
\end{proposition}
\begin{proof}
In the proof of \cite[Proposition~2.1.4]{PartI}, take a common
semipositive twist $\bar U=f^*\bar N$ and put $M=L+U$. Each positive
mixed measure occurring in the expansion has mass $M^{d-k}U^k(f^*B)^e$.
All these masses vanish because $v_B(L)=0$,
$L^{d-1}U(f^*B)^e=L_\eta^{d-1}(N\cdot B^e)=0$, and
$U^2(f^*B)^e=0$ by dimension. For the mass formula see
\cite[Appendix~1]{YZ17} and \cite[\S1.2.8]{CL10}.

In the smooth case, first give $B$ a smooth positive curvature form
$\omega$. For any of the interpolating metrics let $\alpha$ denote
its curvature. Wedging $\alpha^d$ with $f^*\omega^e$ sums the
determinants on the vertical tangent space enlarged by one lifted
$\omega$-orthonormal base direction at a time, up to a positive
factor. The vertical block is semipositive and singular by
$L_\eta^{d-1}=0$. Each determinant is therefore nonpositive by the
matrix calculation in \cite[proof of Proposition~2.1.4]{PartI}.
Since $v_B(L)=0$, the resulting measure has mass zero.
Uniform approximation of the metric on $B$ by smooth positive
metrics and continuity of mixed measures give $\alpha^d\beta^e=0$.
Polarizing along the segment between $h_0$ and $h_1$ gives the result.
\end{proof}

\subsection{Classical functionals}\label{ss:metric-examples}

Let $f$ be smooth over $\BC$ and let $\omega=c_1(\bar B)$ be
K\"ahler. We first compare the functional of Feng--Liu--Wan with
\eqref{eq:H-energy} for $f$-ample $L$. When $v_B(L)=0$, the only
change in \cite[(2.2.1)]{PartI} is the factor $\omega^e/e!$ on $Y$.
For metrics of strictly positive vertical curvature with local weights
$\phi_i$ we obtain
\begin{equation}\label{eq:FLW-normalization}
 \mathcal L_{\rm FLW}(\phi_1,\phi_0)
 =\frac{(2\pi)^d}{d(d+1)e!}\,\Psi_{A,B}(h_0,h_1).
\end{equation}
Indeed $\sqrt{-1}\partial\bar\partial\phi_i=2\pi\alpha_i$ in
\cite[(1.10), (1.13)--(1.14)]{FLW}. Thus
$\alpha^d\beta^e=0$ is the geodesic--Einstein equation at slope
zero, and each solution in $\CH^+(L)$ minimizes the energy there
\cite[Proposition~1.4 and Theorem~0.1]{FLW}. For example, if
$mL=g^*J$ with $J$ ample on a variety of dimension $d-1$, curvature
forms pulled back from $J$ have all products of $d$ factors zero.
Any two such metrics have zero energy, and the corresponding adelic
metrics on compatible models have zero inverse Deligne height. A
pulled-back Fubini--Study metric of strictly positive vertical
curvature is a minimizer, also among smooth fiberwise-semipositive
metrics by interpolation.

For arbitrary slope put $t=v_B(L)/(dL^{d-1}(f^*B)^{e+1})$ and
$(L',h_i')=(L,h_i)-tf^*\bar B$. Relative ampleness gives a positive
denominator. Since $\beta^{e+2}=0$, we have $v_B(L')=0$ and
\[
 (\alpha-t\beta)^d\beta^e
   =\alpha^d\beta^e-dt\,\alpha^{d-1}\beta^{e+1}.
\]
The slope in \cite[(1.15)]{FLW} is $2\pi(e+1)t$. Thus the two first
variations agree with the factor in \eqref{eq:FLW-normalization}, and
both functionals vanish on the diagonal. This proves
\eqref{eq:FLW-normalization} for the original $\phi_0,\phi_1$ with
$\Psi_{A,B}(h_0',h_1')$ on the right.

Now let $X=\BP_Y(E)$ be the bundle of quotient lines of a bundle of
rank $r=d$, put $\xi=\CO_X(1)$ and $L=\xi-r^{-1}f^*\det E$. For a
smooth Hermitian squared metric $G$ on $E$, let $q_G$ denote the
induced metric on $L$. Along a path $G_t$ put $q_t=q_{G_t}$ and
\[
 V_t=G_t^{-1}\dot G_t,\qquad
 R_t=\frac{\sqrt{-1}}{2\pi}\bar\partial(G_t^{-1}\partial G_t),
 \qquad T^0=T-\frac{\operatorname{Tr}T}{r}\operatorname{Id}.
\]
Since $f_*c_1(L)^r=0$, the energy involves only $\bar B$. The
quotient and determinant factors subtract the traces in both the
variation and the horizontal curvature. Thus the fiber integral
computed in \cite[proof of (2.2.2)]{PartI} with $V_t^0,R_t^0$
gives the $(1,1)$-form $\operatorname{Tr}(V_t^0R_t^0)$ on $Y$.
Multiplying by $\omega^e$ and integrating yields
\begin{equation}\label{eq:normalized-bundle-donaldson}
 \Psi_{A,B}(q_0,q_1)=
 \int_0^1\!\int_Y\operatorname{Tr}(V_t^0R_t^0)\wedge\omega^e\,dt .
\end{equation}
At fixed determinant this is Donaldson's functional
\cite{Do85,Do87,Mo94}. Its critical equation is
$\Lambda_\omega R^0=0$. With determinant of constant mean curvature
this is the Hermitian--Einstein equation. Such a metric exists
precisely when $E$ is polystable \cite{NS65,Do83,Do85,UY86}.
Scalar multiplication of
$G$ by a positive function preserves $q_G$. If $\det G_0$ has
constant mean curvature, write $u=\log(\det G/\det G_0)$ and
$\widetilde G=e^{-u/r}G$. Integration by parts on the scalar segment
$e^{tu/r}\widetilde G$ shows that the excess in Donaldson's functional
over its value at $\widetilde G$ is a positive multiple of
$\int_Y\sqrt{-1}\partial u\wedge\bar\partial u\wedge\omega^e$.

\begin{proposition}\label{prop:semistable-projective-infimum}
Suppose that $E$ is slope semistable with respect to $B$ on the smooth
complex projective variety $Y$ and that $\omega=c_1(\bar B)$ is
K\"ahler. For every integrable reference metric $h_0$ on
$L=\CO_{\BP_Y(E)}(1)-r^{-1}f^*\det E$, we have
\[
 \inf_{h\in\BH^+(L)}\Psi_{A,B}(h_0,h)
   =\inf_G \Psi_{A,B}(h_0,q_G)>-\infty ,
\]
where $G$ ranges over the smooth Hermitian metrics on $E$.
\end{proposition}
\begin{proof}
By the cocycle identity we may take $h_0=q_{G_0}$, with $\det G_0$
of constant mean curvature. Let $G_t$ be the solution of Donaldson's
heat equation starting at $G_0$. The flow preserves $\det G_0$.
We write $q_t=q_{G_t}$, $K_t=\Lambda_\omega R_t^0$ and
$a(t)=\|K_t\|_{L^2(Y)}$, using $G_t$ on $\operatorname{End}(E)$ and
the volume form of $\omega$ on $Y$. By
\eqref{eq:normalized-bundle-donaldson}, $\Psi_{A,B}(q_0,q_t)$ is Donaldson's
functional along the flow. Semistability bounds it below, and its
derivative is a negative constant times $a(t)^2$
\cite[Theorem~6.8.25, Proposition~6.9.1 and Theorem~6.10.13]{Kob87}.
Thus $\Psi_{A,B}(q_0,q_t)$ decreases to a finite limit $D_\infty$, and
$\int_0^\infty a(t)^2\,dt<\infty$. We choose $t_j\to\infty$ with
$\sqrt{t_j}\,a(t_j)\to0$.

We compare the quotient metrics $q_t$ with an arbitrary
$h\in\CH^+(L)$. Let $\phi_t$ and $\phi$ be the local weights of
$q_t$ and $h$, and put $b_t(y)=\sup_{X_y}|\phi-\phi_t|$. We claim that
\begin{equation}\label{eq:quotient-energy-comparison}
 \Psi_{A,B}(q_t,h)\ge -C\,a(t)\,\|b_t\|_{L^2(Y)}
\end{equation}
with $C$ independent of $t$ and $h$. Since $f_*c_1(L)^r=0$, the slope
is zero and \eqref{eq:FLW-normalization} identifies the two functionals
up to a positive factor. For the rational bundle $L$, we apply the
formula to $rL$ and rescale the weights by $1/r$.

Fix $t$ and join $\phi_t$ to $\phi$ by their approximate geodesic
$\Phi_{s,\varepsilon}$, $0\le s\le1$, using the fiber volume form of
$q_t$ in the approximation. On each $X_y$, put
$m_y=\min_{X_y}(\phi-\phi_t)$. The function
$\phi_t+s m_y+c\varepsilon s(s-1)$ is a lower barrier for the
fiber equation when $c$ is sufficiently large. Comparison at $s=0$
and convexity in $s$ give
\[
 m_y-c\varepsilon
 \le \partial_s\Phi_{0,\varepsilon}
 \le \phi-\phi_t\qquad\text{on }X_y .
\]
The estimate of Feng--Liu--Wan bounds $\Psi_{A,B}(q_t,h)$ below by its
initial derivative along this path up to an error tending to zero
with $\varepsilon$ \cite[Lemma~1.5 and (1.23)]{FLW}.
At a quotient line over $y$ with orthogonal projector $\Pi$, the
$\omega$-trace of the horizontal curvature of $q_t$ is
$\operatorname{Tr}(\Pi K_t(y))$. Its absolute value is at most
$|K_t(y)|$. In the first variation the vertical volume form has
mass one on $X_y$. Letting $\varepsilon$ tend to zero with the
endpoints fixed therefore gives
$\Psi_{A,B}(q_t,h)\ge -C\int_Y b_t|K_t|\,\omega^{e+1}$.
Cauchy--Schwarz gives \eqref{eq:quotient-energy-comparison}.

The heat equation also gives
$\sup_{X_y}|\dot\phi_t|\le C|K_t(y)|$. Integrating in time, we obtain
\[
 \|b_t\|_{L^2(Y)}
 \le \|b_0\|_{L^2(Y)}+C\int_0^t a(s)\,ds
 \le \|b_0\|_{L^2(Y)}
      +C\sqrt t\left(\int_0^\infty a(s)^2\,ds\right)^{1/2}.
\]
It follows from \eqref{eq:quotient-energy-comparison} that
$\liminf_j \Psi_{A,B}(q_{t_j},h)\ge0$. The cocycle identity now gives
$\Psi_{A,B}(q_0,h)\ge D_\infty$ for every $h\in\CH^+(L)$.
Lemma~\ref{lem:relative-nef-approximation} extends this bound to
$\BH^+(L)$. Since the quotient metrics $q_t$ themselves belong to
$\CH^+(L)$ and $\Psi_{A,B}(q_0,q_t)\to D_\infty$, both infima equal $D_\infty$.
\end{proof}

For the unnormalized bundle $\xi$ with $c_1(E)\cdot B^e=0$, the
same moment identity, now for $V_t,R_t$, gives the mixed version of
\cite[(2.2.2)]{PartI}, namely
\begin{equation}\label{eq:Donaldson-normalization}
 \Psi_{A,B}(h_0^\xi,h_1^\xi)
   =\mathcal D_0(G_0,G_1)+\frac{r+1}{r}\mathcal D_{\det}(G_0,G_1),
\end{equation}
where $h_i^\xi$ are the quotient metrics and $\mathcal D_0$ is the
right side of \eqref{eq:normalized-bundle-donaldson}. Writing
$u=\log(\det G_1/\det G_0)$, we have
\[
 \mathcal D_{\det}
 =\int_0^1\!\int_Y\operatorname{Tr}(V_t)\operatorname{Tr}(R_t)
                       \wedge\omega^e\,dt
 =\tfrac12\int_Y u\bigl(c_1(\det E,\det G_0)+c_1(\det E,\det G_1)\bigr)
                       \wedge\omega^e .
\]
In particular the determinant term vanishes at fixed determinant. For
the relation to approximate geodesic--Einstein structures under long-time
smooth existence of the flow, see \cite[Theorem~2.8]{Wan}.

For nef, $f$-ample $L$ of mixed degree $v_B(L)=0$, the approximate
critical metrics of \cite[\S2.2]{PartI} have the following estimate.
We choose a K\"ahler form $\Omega_\epsilon$ in the ample class
$c_1(L)+\epsilon[\beta]$ and a metric on $L$ with curvature
$\alpha_\epsilon=\Omega_\epsilon-\epsilon\beta$. Its vertical
curvature is positive. Since $\beta^{e+2}=0$,
$\alpha_\epsilon^d\beta^e$ is the difference of the two positive
measures $\Omega_\epsilon^d\beta^e$ and
$d\epsilon\Omega_\epsilon^{d-1}\beta^{e+1}$, both of mass
$d\epsilon L^{d-1}(f^*B)^{e+1}$. Hence
$\int_X|\alpha_\epsilon^d\beta^e|\le2d\epsilon L^{d-1}(f^*B)^{e+1}\to0$.
Finski computes the infimum of this integral more generally using
asymptotic slopes \cite[Theorem~1.1]{Finski26WZW}. See also
\cite{Finski26YM,Finski26HN}.

\subsection{Arithmetic lower bounds}\label{ss:energy}

For $f:X\to Y$ as in Theorem~\ref{thm:main-height}, suppose that
$L|_{X_C}$ is nef. Recall that
\begin{equation}\label{eq:H-arithmetic-deficit}
 I_{A,B}(\bar L)=\widehat{\deg}\,D_{A,B}(\bar L)
 =(d+1)v_B(L)(\bar L\cdot\bar A^d\cdot(f^*\bar B)^e)
          -\bar L^{d+1}\cdot(f^*\bar B)^e .
\end{equation}
Taking the arithmetic degree of the quotient of the two metrics on
$D_{A,B}(L)$ gives
\begin{equation}\label{eq:H-height-energy-variation}
 I_{A,B}(\bar L_1)-I_{A,B}(\bar L_0)
       =\sum_v\Psi_v(h_{0,v},h_{1,v}),
\end{equation}
with the place weights included as in \cite[(2.3.2)]{PartI}.
Write $\bar L_0^v(h_v)$ for the adelic metric equal to $h_v$ at $v$
and to $\bar L_0$ elsewhere.

\begin{thm}\label{thm:lb}
If $\bar L_0^v(h_v)$ is integrable and relatively nef, then
$\Psi_v(h_{0,v},h_v)\ge -I_{A,B}(\bar L_0)$.
\end{thm}
\begin{proof}
Theorem~\ref{thm:main-height} gives $I_{A,B}(\bar L_0^v(h_v))\ge0$,
and \eqref{eq:H-height-energy-variation} has a single term.
\end{proof}

\begin{lemma}\label{lem:one-place-globalization}
Let $f:X\to Y$ be any projective morphism of projective varieties
over $K$ and let $L$ be $f$-ample. If $\bar L_0$ is strongly
relatively nef, so is $\bar L_0^v(h_v)$ for every strongly relatively
nef metric $h_v$ on $L_v$.
\end{lemma}
\begin{proof}
As in the proof of \cite[Lemma~2.3.2]{PartI}, we use ampleness of the base
twist on the projective variety $Y$. A sufficiently large multiple of
that twist dominates a semipositive ample twist for $h_v$ and gives
an ample underlying bundle on $X$. The common constant over a number
field, or the common divisor on $T$ over $k(T)$, makes both adelic
factors nef and cancels in their difference. This proves the assertion
over $Y$ of any dimension.
\end{proof}

\begin{cor}\label{cor:H-local-bound}
Under the hypotheses of Theorem~\ref{thm:main-height}, suppose that
$L$ is $f$-ample and $\bar L_0$ is strongly relatively nef. At every
place $v$ we have
\[
 \inf_{h\in\BH^+(L_v)}\Psi_v(h_{0,v},h)
            \ge-I_{A,B}(\bar L_0)>-\infty .
\]
\end{cor}
\begin{proof}
Lemma~\ref{lem:one-place-globalization} and Theorem~\ref{thm:lb}
give the bound for strongly relatively nef metrics. We apply
Lemma~\ref{lem:relative-nef-approximation} to pass to $\BH^+(L_v)$.
At a real place, taking the twists together with their conjugates
and averaging the approximating weights of their semipositive factors
gives conjugation-invariant approximants.
\end{proof}

\section{Families of fibrations}\label{sec:families}\label{ss:family-energy}
We prove Theorem~\ref{thm:main-energy} by applying the height inequality
to the fibers of a family $X\to Y\to S$. When $\dim Y_s=e+1$, the
Deligne pairing on $S$ contains $e$ copies of the polarization from
$Y$. To apply the height inequality, we find algebraic parameters $s$
for which the perturbed bundle is ample on the inverse image of a
general curve section of $Y_s$.

\subsection{Models, polarizations, and specialization}
\label{ss:adelic-models-approximations}\label{ss:family-critical}

Let $K$ be a global field, and let $T=\Spec\mathcal O_K$ for a number
field or the smooth projective curve of $K$ over its finite constant
field otherwise. Over a smooth geometrically integral quasi-projective
$K$-variety $S$, we consider a projective tower
\begin{equation}\label{eq:tower}
 X\xrightarrow{f}Y\xrightarrow{\pi}S,\qquad g=\pi\circ f .
\end{equation}
We assume that $f$ is flat and surjective of relative dimension $d-1$,
and that $\pi$ and $g$ are smooth and surjective with geometrically
integral fibers of dimensions $e+1$ and $d+e$, where $e\ge1$.
Let $\xi$ be the generic point of $S$. We assume the generic fiber
of $f_\xi$ is geometrically integral, and fix an $f$-ample rational
line bundle $L$.

\begin{lem}\label{lem:polarized-arithmetic-spread}
We may extend \eqref{eq:tower} to a flat projective tower of integral schemes
$\mathcal X\to\mathcal Y\to\mathcal S\to T$, with $\mathcal S$ normal
and $S^c=\mathcal S_K$ a projective compactification of $S$, smooth
in characteristic zero. Moreover, we may extend $L$ to a relatively ample
rational bundle $\mathcal L$ on $\mathcal X/\mathcal Y$.
The extension $\mathcal L$ gives a strongly relatively nef adelic reference $\bar L_0$
with model metrics at the finite places and smooth, vertically
positive metrics in the ambient sense at infinity. For
$s\in S(K')$, with $K'/K$ finite, the pullback to the normalization
of $T$ in $K'$ gives integral flat models of $X_s\to Y_s$ with
$\mathcal L_s$ relatively ample.
\end{lem}
\begin{proof}
The construction in \cite[Lemma~3.2.1]{PartI} embeds $X$ over $Y$
by a multiple of $L$, takes the projective closure, and flattens
first over the model of $Y$ and then over the model of $S$. All
three operations apply in every relative dimension of $Y/S$.
On the resulting tower, an ample twist from $\mathcal Y$ gives
the reference metric as in that proof. Flatness and geometric
integrality of $X_s,Y_s$ give the last assertion.
\end{proof}

For Theorem~\ref{thm:main-energy}, take $K$ to be the algebraic
closure of $K_0$ in $F$ and choose the tower over $S$ to spread out
$X_F\to Y_F$, so that $K(S)=F$. Fix nef rational model bundles
$\bar A_0$ on the total variety and $\bar B$ on the intermediate
variety, with ample underlying bundles and, at infinity, smooth
semipositive metrics compatible with conjugation.

\begin{defn}\label{def:admissible-polarizations}
We call $(\bar A_\sigma,\bar B_\sigma)_{\sigma:F\hookrightarrow E}$
\emph{admissible} if there is one polarized model $(\bar A_0,\bar B)$ as above
induces the whole family. For each $K_0$-embedding $\sigma$ we restrict
the model bundles to the fibers over the corresponding $K$-generic
point of $S_\iota(E)$, where $\iota=\sigma|_K$. We take $\bar B_\sigma$
to be the restriction of $\bar B$ and normalize $\bar A_{0,\sigma}$ by
\[
 \bar A_\sigma
 =\bigl(A_{0,\sigma}^d\cdot(f_\sigma^*B_\sigma)^e\bigr)^{-1/d}
      \bar A_{0,\sigma}.
\]
\end{defn}

To realize the fixed polarizations and $\mathcal L$ on one tower,
we use the common model constructed after
\cite[Definition~3.2.2]{PartI}. Its total variety dominates the
prescribed model of $X$ and its intermediate variety dominates the
prescribed model of $Y$. Thus we can pull back both $\bar A_0$ and
$\bar B$ while using the relatively ample extension of $L$
constructed on the graph there.

For $s\in S$ put $v_{B,s}(M)=M_s^d\cdot(f_s^*B_s)^e$ and
$\bar A=v_{B,s}(A_0)^{-1/d}\bar A_0$. Flatness makes the
coefficient independent of $s$. On $S^c$ we give the line $D_{A,B}(L)$ of the introduction the
Deligne metric induced by $\bar L$, $\bar A$ and $\bar B$.

The $d$ entries $\bar A$ contribute the rational factor
$1/v_{B,s}(A_0)$. For a finite extension $K'/K$ and $s\in S(K')$, base change gives
\begin{equation}\label{eq:B-mixed-pairing}
 h_{D_{A,B}(\bar L)}(s)=\frac{I_{A_s,B_s}(\bar L_s)}{[K':K_0]} .
\end{equation}
The change of its local
weight is $\Psi_s=\Psi_{A_s,B_s}$.

Fix $a_1,\ldots,a_e$ with $a_iB_\xi$ very ample, and write $C$ for
the complete-intersection curve over the geometric generic point of
$\prod_i|a_iB_\xi|$. Put $X_{\xi,C}=X_\xi\times_{Y_\xi}C$.

\begin{lem}\label{lem:B-specialization}
Suppose that $L|_{X_{\xi,C}}$ is nef. For each rational $t>0$ there
is a dense $K$-open $U_t\subset S$ such that for every geometric
$s\in U_t$, the generic fiber of $f_s$ is geometrically integral
and $L_s+tA_{0,s}$ is ample on the inverse image of the geometric
generic curve cut out by $|a_1B_s|,\ldots,|a_eB_s|$.
\end{lem}
\begin{proof}
After shrinking $S$, the $a_iB$ are relatively very ample and
$\pi_*(a_iB)$ are locally free and commute with base change. The
product over $S$ of the projective bundles
$\BP_S((\pi_*(a_iB))^\vee)$ is a parameter space $W\to S$ with
fiber $W_s=\prod_{i=1}^e |a_iB_s|$.
Over a dense open of $W$ the universal intersection is a smooth
geometrically integral curve in each fiber. On the inverse image in
$X$ of the geometric generic curve, $L+tA_0$ is ample by nefness
of $L$. Openness of ampleness for this projective family gives a
nonempty open $W_t\subset W$ where the restriction is relatively
ample. The image of $W_t$ in $S$ is open because $W\to S$ is smooth.
Let $Y^\circ\subset Y$ be a dense open over which $f$ has geometrically
integral fibers. We take $U_t$ to be the intersection of the image of
$W_t$ with $\pi(Y^\circ)$. It is dense and open since $\pi$ is smooth.
For every geometric $s\in U_t$, the open $Y^\circ_s$ contains the generic
point of $Y_s$, and $(W_t)_s$ contains the geometric generic point of $W_s$.
\end{proof}

\begin{lem}\label{lem:metric-extension}\label{lem:nonarch-metric-extension}
Fix an embedding $\iota:K\hookrightarrow\BC$ and write
$S_\iota=S\times_{K,\iota}\BC$. For every $s\in S_\iota(\BC)$,
each $h_s\in\CH^+(L_s)$ extends to a smooth vertically positive
metric on $L_\iota$ equal to a given smooth vertically positive
reference $h_0$ outside $g_\iota^{-1}(K_s)$ for some compact
$K_s\subset S_\iota(\BC)$. Over a nonarchimedean local field $E$, for a
flat projective $f:X\to Y$ in a tower of projective $E$-varieties
and $s\in S(E)$ with $Y_s$ smooth geometrically integral and $X_s$
geometrically integral, every $h_s\in\CH^+(L_s)$ extends to a
relatively ample rational model metric on a flat model of $f$.
In either case the restrictions of these extensions approximate
every $h_s\in\BH^+(L_s)$ uniformly and in $\Psi_s(h_{\rm ref},-)$
for every integrable reference metric $h_{\rm ref}$.
\end{lem}
\begin{proof}
At a complex place we use the cutoff in \cite[Lemma~3.2.3]{PartI}.
It is pulled back from $S$, so its derivatives vanish on
$\ker df\subset\ker dg$. At a real place the extension respects
conjugation by averaging the weights or by using disjoint conjugate
neighborhoods.

At a nonarchimedean place we follow the proof of
\cite[Lemma~3.2.4]{PartI}. We need an ample model
of $N_s=N|_{Y_s}$ on a model dominating the base of $h_s$, where
$N$ is ample on $Y$. In higher dimension we obtain it as follows.
Let $\mathcal Y_1$ be the base of a relatively ample model inducing
$h_s$. The embedding by $aN_s$ for $a$ sufficiently divisible gives
a model $(\mathcal Y_0,\mathcal P)$ with $\mathcal P$ ample.
An admissible blowup $p:\mathcal Y'\to\mathcal Y_0$ dominates the
closure of the graph in $\mathcal Y_0\times_{\mathcal O_E}\mathcal Y_1$
\cite[Tag 081Q]{Stacks}. If $\mathcal J$ is the invertible exceptional
ideal, then $\mathcal J\otimes p^*\mathcal P^{\otimes q}$ is ample for
$q$ large and extends $aqN_s$. Dividing by $aq$ gives the required
ample rational model $\mathcal N_s$. On the strict transform over
$\mathcal Y'$ the pullback of the model bundle inducing $h_s$
remains relatively ample. Adding a sufficiently large multiple of
the pullback of $\mathcal N_s$, we obtain $h_s=a_s-f_s^*n_s$, where
$a_s$ and $n_s$ are ample model metrics on $(L+mf^*N)_s$ and $mN_s$,
respectively. The lifting of their sections and the flat graph model
in \cite[Lemma~3.2.4]{PartI} give the extension.

For $h_s\in\BH^+(L_s)$, Lemma~\ref{lem:relative-nef-approximation}
gives uniform approximation by metrics in $\CH^+(L_s)$ with
convergence of the energies for the fixed $e$ factors $f_s^*\bar B_s$.
We extend those metrics by the preceding construction.
\end{proof}

Consequently the restrictions give the same infimum of the energy,
and hence the same supremum of the induced norms on $D_{A,B}(L)_s$,
as the full class $\BH^+(L_s)$. To replace a strongly relatively
nef adelic reference at a single place we use
Lemma~\ref{lem:one-place-globalization}.

\subsubsection*{The envelope on the parameter space}

We formulate the analogue of \cite[Conjecture~3.1.1]{PartI} for
$X\xrightarrow{f}Y\to S$. First consider a tower as above over
$F=\BC$ or a complete discretely valued field. We fix semipositive
metrized polarizations $\bar A,\bar B$ with $v_{B,s}(A)=1$. We use
the family metric classes $\CH^+(L)$ of \cite[\S3.1]{PartI}.
Over $\BC$ we use smooth integrable metrics on $L$ with strictly positive
curvature on $\ker df$, and over a nonarchimedean field we use
strongly relatively nef rational model metrics. For $s\in S^{\an}$
and $u\in D_{A,B}(L)_s$, put
\begin{equation}\label{eq:mixed-family-envelope}
 |u|_{\mathrm{env},s}=\sup_{h\in\CH^+(L)}|u|_{D_{A,B}(h),s}.
\end{equation}
Relative to a fixed $h_0\in\CH^+(L)$, the change of weight from
$D_{A,B}(h_0)$ to this envelope is
$\inf_{h\in\CH^+(L)}\Psi_{A,B}(h_0,h)(s)$.

Over a number field $K$, we fix nef model polarizations as above
and use conjugation-invariant metrics at the real places. In the
height infimum below, $\CH^+(L)$ denotes the strongly relatively nef
adelic metrics on $L$ induced by Hermitian models, with smooth,
strictly positive vertical curvature at infinity.

\begin{conj}\label{conj:mixed-family-envelope}
Suppose that $L$ is $f$-ample. For each geometric point $s$ of $S$,
assume that $f_s$ has geometrically integral generic fiber and that
$L_s$ is nef on $f_s^{-1}(C_s)$ for some geometric generic curve
section $C_s$ of $(Y_s,B_s)$ as in Theorem~\ref{thm:main-height}.
\begin{enumerate}[label=(\roman*),leftmargin=*,itemsep=0.3\baselineskip]
\item For a local family over $F$, the envelope
\eqref{eq:mixed-family-envelope} is finite at every $s\in S^{\an}$
and defines a continuous semipositive metric on $D_{A,B}(L)$.
\item For a family over a number field $K$, the local envelopes
define a nef adelic line bundle $D_{A,B}^{\mathrm{env}}(L)$ on $S$.
For every $s\in S(\bar K)$ its height is
\[
 h_{D_{A,B}^{\mathrm{env}}(L)}(s)
   =\inf_{\bar L\in\CH^+(L)}h_{D_{A,B}(\bar L)}(s).
\]
\end{enumerate}
\end{conj}

\subsubsection*{The average energy at a place}

We use Theorem~\ref{thm:main-height} on algebraic fibers to bound the
average of the local energy at a fixed place $v$. We first express the
average over the conjugates of an algebraic point as a difference of
heights. We then pass to a limiting measure on the parameter space.
The resulting lower bound \eqref{eq:average-mixed-energy} is independent
of the metric on $L$.
Let $X^c\to Y^c\to S^c$ be the compactified tower of
Lemma~\ref{lem:polarized-arithmetic-spread}, with the prescribed
polarizations $\bar A_0,\bar B$ and reference bundle $\bar L_0$.
Suppose that $L$ is nef on $X_{\xi,C}$ as in
Lemma~\ref{lem:B-specialization}. Fix a place $v$ of $K$ and write
$W=(S^c)_{K_v}^{\an}$. For a complex embedding we take
$W=S^c(\BC)$, passing first to $K(i)$ if the embedding of $K$ is real.
We denote the resulting field again by $K$. The chosen complex
parameter space and its polarizations are preserved by this base change.

Fix a nonnegative ample height $h_D$ on $S^c$. At a complex place,
\cite[Lemma~3.3.2]{PartI} gives a probability measure $\mu$ on $W$
dominating a positive multiple of any prescribed smooth positive volume
form on $S^c(\BC)$. At a nonarchimedean place,
\cite[Lemma~3.4.3]{PartI} gives a probability measure on $W$ whose
restriction to $S(K_v)$ dominates the Haar measure class. In both cases
there is $H>0$ such that, for every dense $K$-open $U\subset S$, $\mu$
is a weak limit of convex combinations of the normalized orbit measures
$\nu_{a,v}$ for $a\in U(\bar K)$ with $h_D(a)\le H$. For $K'=K(a)$
we use $\nu_{a,v}=[K':K]^{-1}\sum_{w\mid v}[K'_w:K_v]\delta_{a_w}$,
where $a_w$ is the local point induced by $w$. At a complex place
we use the chosen embedding of $K$; at a nonarchimedean place we
take $K'/K$ separable. The limiting measure gives mass zero to $S^c\setminus S$.

At a complex place we consider smooth metrics $h$ with strictly
positive vertical curvature that agree with $h_0$ off $g^{-1}(K_h)$
for some compact $K_h\subset S(\BC)$. At a nonarchimedean place we
take all strongly relatively nef metrics $h$ on $L$ over $X^c_{K_v}/Y^c_{K_v}$.
We write $E_h=\Psi_{A,B}(h_0,h)$ on $W$. At a complex place we
extend the change of metric by zero near the boundary and make the
conjugate change at the conjugate embedding. A sufficiently positive
base twist makes this extended metric semipositive, as in
\cite[proof of Lemma~3.3.4]{PartI}. In both cases Lemma~\ref{lem:one-place-globalization} allows us to replace the
$v$-component of $\bar L_0$ by $h$.

For rational $0<\epsilon\le1$, we put
\[
 \bar L_\epsilon(h)=\bar L_0^v(h)+\epsilon\bar A_0,\qquad
 \bar R_{\epsilon,h}=D_{A,B}(\bar L_\epsilon(h)),\qquad
 v_\epsilon=(L_\xi+\epsilon A_{0,\xi})^d\cdot(f_\xi^*B_\xi)^e .
\]
The coefficient $v_\epsilon$ is the mixed degree on every parameter fiber by
flatness. Lemma~\ref{lem:B-specialization} gives a dense open $U_\epsilon$
on which $L+\epsilon A_0$ is ample on the inverse image of the geometric
generic curve section of each $Y_a$. Thus Theorem~\ref{thm:main-height}
and \eqref{eq:B-mixed-pairing} give $h_{\bar R_{\epsilon,h}}(a)\ge0$ for
$a\in U_\epsilon(\bar K)$. By multilinearity, $\bar R_{\epsilon,h_0}$ is a polynomial
in $\epsilon$ with finitely many fixed adelic line bundles as coefficients.
Ample height comparison gives
$h_{\bar R_{\epsilon,h_0}}(a)\le c(1+h_D(a))$ uniformly in $0<\epsilon\le1$.

Let $E_{\epsilon,h}$ be the weight of $\bar R_{\epsilon,h}$ minus that
of $\bar R_{\epsilon,h_0}$ at $v$. These metrics agree away from $v$.
For $a\in U_\epsilon(\bar K)$ put $K'=K(a)$. By the normalization in
\eqref{eq:B-mixed-pairing}, the difference of heights is
\[
 \begin{aligned}
 h_{\bar R_{\epsilon,h}}(a)-h_{\bar R_{\epsilon,h_0}}(a)
 &=\frac{1}{[K':K_0]}\sum_{w\mid v}[K'_w:K_v]\,E_{\epsilon,h}(a_w)\\
 &=\frac{1}{[K:K_0]}\int_W E_{\epsilon,h}\,d\nu_{a,v}.
 \end{aligned}
\]
The first equality sums the local changes of weight at the places of
$K'$ above $v$. The second follows from the definition of $\nu_{a,v}$.
For $h_D(a)\le H$ the height estimates above bound this difference
below by $-c(1+H)$. For fixed $\epsilon$ the
orbit measures from $U_\epsilon$ approximate $\mu$ in the sense above, so
the same bound holds after integration against $\mu$. For fixed $h$,
the functions $E_{\epsilon,h}$ are polynomials in $\epsilon$ with continuous
coefficients on $W$. These coefficients are differences of weights of the mixed
Deligne pairings defining $\bar R_{\epsilon,h}$, which have the same underlying
lines for $h$ and $h_0$. Their constant term is $E_h$. Thus
$E_{\epsilon,h}\to E_h$ uniformly on $W$. We have proved the common bound
\begin{equation}\label{eq:average-mixed-energy}
 \int_W E_h\,d\mu\ge-b,\qquad b=[K:K_0]c(1+H).
\end{equation}
The restriction to $U_\epsilon$ precedes passage to $\mu$ for each $\epsilon$.
In particular, the bound is uniform in $h$ although $U_\epsilon$ varies with $\epsilon$.

\subsection{Extension to complex points}\label{ss:complex-points}

Let $K$ be a number field and fix $K\hookrightarrow\BC$. For the family
and polarizations of Section~\ref{ss:adelic-models-approximations}, fix a
smooth metric $h_0$ on $L$ with strictly positive vertical curvature and put
\[
 m_\infty(s)=\inf_{h_s\in\CH^+(L_s)}\Psi_s(h_{0,s},h_s),
 \qquad s\in S(\BC).
\]

\begin{thm}\label{thm:ae-energy}
If $L|_{X_{\xi,C}}$ is nef, then $m_\infty$ is upper semicontinuous
and belongs to $L^1_{\rm loc}(S(\BC))$. The infimum is the same
over $\BH^+(L_s)$, and, when $f$ is smooth, over all smooth
fiberwise-semipositive metrics. For the model reference of
Lemma~\ref{lem:polarized-arithmetic-spread}, $m_\infty$ is integrable
for every smooth positive volume form on the compactification $S^c(\BC)$,
restricted to $S(\BC)$.
\end{thm}

\begin{proof}
We first take the model reference. By Lemma~\ref{lem:metric-extension},
$m_\infty$ is the infimum of the $E_h$ in \eqref{eq:average-mixed-energy}
at each $s$. The approximation assertion in that lemma gives the
infimum over $\BH^+(L_s)$. For smooth $f$, interpolation with $h_{0,s}$
approximates each smooth fiberwise-semipositive metric by metrics with
strictly positive vertical curvature.

We apply the cutoff argument of \cite[Lemma~3.3.3]{PartI} to $E_h$.
The cutoff functions are pulled back from $S$, so the glued metric
has positive curvature on $\ker df$. The additional factors
$c_1(f_s^*\bar B_s)^e$ are fixed smooth forms. Hence the energies of
the finitely many interpolating paths are continuous and uniformly
bounded on the common compact parameter support. The finite minima of
the $E_h$, including $E_{h_0}=0$, can therefore be approximated in
$L^1(\mu)$ by single test energies as in that lemma. See also
\cite[Theorem~14.60]{RW98}.
Its countable-subcover argument and \eqref{eq:average-mixed-energy}
give $\int(-m_\infty)\,d\mu\le b$. Domination of smooth volume gives
the asserted integrability. The infimum of the continuous $E_h$ is
upper semicontinuous. Changing the smooth reference to $h'_0$ adds the
smooth function $s\mapsto\Psi_s(h'_{0,s},h_{0,s})$. This preserves
local integrability, and global integrability when $S$ is projective.
\end{proof}

\begin{thm}\label{thm:intrinsic-energy}
For $K_0=\BQ$ and $E=\BC$, the bound of
Theorem~\ref{thm:main-energy} holds outside one set of measure zero
for all integrable reference metrics. The infimum equals the one over
$\CH^+(L_\sigma)$, and, for smooth $f$, the one over all smooth
fiberwise-semipositive metrics.
\end{thm}
\begin{proof}
For each embedding of the constant field $K$ in $\BC$,
Theorem~\ref{thm:ae-energy} gives a finite infimum at almost every
complex point of the fixed spread. Its $K$-generic complex points
are the embeddings of $F$ extending that embedding of $K$ and have
full measure \cite[proof of Theorem~3.3.5]{PartI}. At each such point
the cocycle identity changes the infimum by a finite constant when
we change the integrable reference metric.
\end{proof}
In particular the envelope norm \eqref{eq:mixed-family-envelope} on
$D_{A,B}(L)$ is finite at almost every embedding.

\subsection{Nonarchimedean realizations}\label{ss:nonarchimedean}

Let $K_0=\BQ$ or $\mathbb F_p(t)$ and let $E$ be a finite extension
of $\BQ_p$ or $\mathbb F_p((t))$, respectively, with its fixed
embedding of $K_0$. We use Haar measure in local coordinates on
$S_\iota(E)$ for each embedding $\iota:K\hookrightarrow E$ of the
constant field. Its $K$-generic points have full measure, and the
measure class on embeddings is independent of the model of $F$, as in
\cite[\S3.4]{PartI}.

\begin{thm}\label{thm:padic-intrinsic-energy}
Under the hypotheses of Theorem~\ref{thm:main-energy}, for almost
every $K_0$-embedding $\sigma:F\hookrightarrow E$ and every integrable
reference metric $h_{0,\sigma}$, we have
\[
 \inf_{h\in\BH^+(L_\sigma)}\Psi_\sigma(h_{0,\sigma},h)
 =\inf_{h\in\CH^+(L_\sigma)}\Psi_\sigma(h_{0,\sigma},h)>-\infty.
\]
For the fixed polarized model, the exceptional set is independent
of the reference metric.
\end{thm}
\begin{proof}
Fix $\iota:K\hookrightarrow E$. After a finite extension of constants
inside $E$ we may assume $K_v=E$, as in the proof of
\cite[Theorem~3.4.2]{PartI}. This reduction also treats inseparable
local extensions in positive characteristic and preserves $S_\iota(E)$.
We apply \eqref{eq:average-mixed-energy} on $W=(S^c)_E^{\an}$. We show
that the mixed energy satisfies the estimate needed in the proof of
\cite[Lemma~3.4.1]{PartI} to pass from $E_h$ to its infimum.

For finitely many strongly relatively nef metrics $h_0,\ldots,h_k$,
choose a common semipositive base twist on $N$ over $Y^c_E$ and put
$U=f^*N$, $M=L+U$. The maximum of weights constructed in that proof
defines a strongly relatively nef metric $h_\delta$. On fibers over
a set of $\mu$-measure at least $1-\delta$, it differs by a constant
from a metric realizing $\min_{0\le j\le k} E_{h_j}$ there. The signed
measure in \eqref{eq:H-energy} has mass zero by $v_{B,s}(A)=1$, so
$E_{h_\delta}$ equals the minimum on this set.

On every fiber the selected metric has the same semipositive twist
on $M_s$, and its weight differs from that of $h_{0,s}$ by at most
a constant $c_*$ independent of $\delta$. Expanding the $d$ curvature
factors in \eqref{eq:H-energy} using $L_s=M_s-U_s$ and wedging with
the positive factors $c_1(f_s^*\bar B_s)^e$ gives
\[
 |E_{h_\delta}(s)|\le(d+1)c_*
 \left((M_s+U_s)^d\cdot(f_s^*B_s)^e+|v_{B,s}(L)|\right).
\]
For $s\in W$, let $x$ be its underlying scheme point and let
$\mathcal H(s)$ be its completed residue field. The two degrees in the
bound are $P_s^d\cdot(f_s^*B_s)^e$ for $P=M+U$ and $P=L$. We compute
them on $X_s=X_x\times_{\kappa(x)}\mathcal H(s)$ using the cycle
\[
 [X_s]_{d+e}=\sum_Z
 \operatorname{length}_{\mathcal O_{X_s,\eta_Z}}
      (\mathcal O_{X_s,\eta_Z})[Z],
\]
where $Z$ runs over the reduced irreducible components of dimension
$d+e$ and $\eta_Z$ is the generic point of $Z$. This also specifies
the multiplicity of each component over the boundary $S^c\setminus S$.
Extension of the residue field preserves both degrees. On the algebraic
fibers both degrees are constant on each connected component of the
base. Indeed, after clearing denominators, each mixed degree is $d!e!$
times the coefficient of $a^d b^e$ in
$\chi(X_x,aP_x+b f_x^*B_x)$, and flatness makes this polynomial locally
constant in $x$ \cite[Tags 02QS and 0B9T]{Stacks}.

We may therefore choose a constant bounding $|E_{h_\delta}|$ and
$\big|\min_j E_{h_j}\big|$ on $W$, independently of $\delta$.
Since $E_{h_\delta}=\min_j E_{h_j}$ outside a set of $\mu$-measure
at most $\delta$, we obtain
$\|E_{h_\delta}-\min_j E_{h_j}\|_{L^1(\mu)}=O(\delta)$.

The space $W$ is second countable \cite[Theorem~1.1]{HLP14}.
The countable selection and monotone convergence argument of
\cite[Lemma~3.4.1]{PartI} now gives
$\int_W(-\inf_h E_h)\,d\mu\le b$. By
Lemma~\ref{lem:nonarch-metric-extension}, at $s\in S(E)$ this infimum
equals the infimum in the theorem for the model reference: local
metrics in $\CH^+(L_s)$ extend and their energies approximate all those
in $\BH^+(L_s)$. Domination of Haar measure on $S(E)$ gives the bound
almost everywhere. As at complex places, the cocycle identity gives
the same exceptional set for every integrable reference metric.
\end{proof}

\clearpage
\phantomsection

%------------------------------------------------------
\end{document}